\documentclass[12pt]{article}
\usepackage{amsfonts}
\usepackage{latexsym, times}
\usepackage{graphics}
\usepackage{graphicx,float,latexsym,shortvrb}

\date{}
\begin{document}
\vskip 3mm

\noindent \vskip 3mm

\title{Two-Sided Power Random Variables}
\author{ H. Homei }
\maketitle

\noindent Department of Statistics, Faculty of Mathematical
Sciences, University of Tabriz, P.O. Box 5166617766, Tabriz, Iran.

[Email: homei@tabrizu.ac.ir],
[Phone: 0098 411 339 2901], [Fax: 0098 411 334 2101].\\

\noindent Running title:
\title{Two-Sided Power Random Variables}\\

\noindent Keywords and Phrases: Two-sided power, moment, triangular distribution, power distribution.\\

\noindent Mathematics Subject Classification [2010]: 62E10

\pagebreak

\begin{abstract}
We study a well-known problem concerning a random variable $Z$
uniformly distributed between two independent random variables. Two
different extensions, conditionally directed power distribution and
conditionally undirected power distribution, have been introduced
for this problem. For the second method, two-sided power random
variables have been defined.
\end{abstract}

\section{Introduction}
Van Assche (1987) introduced the notion of a random variable $Z$
uniformly distributed between two independent random variables $X_1$
and $X_2$, which arose in studying the distribution of products of
random $2\times2$ matrices  for stochastic search of global maxima.
By letting $X_1$ and $X_2$ to have identical distributions, he
derived that: (i) for $X_1$ and $X_2$ on $[-1,1]$, $Z$ is uniform on
$[-1,1]$ if and only if $X_1$ and $X_2$ have an Arcsin distribution;
and (ii) $Z$ possesses the same distribution as $X_1$ and $X_2$ if
and only if $X_1$ and $X_2$ are degenerated or have a Cauchy
distribution. Soltani and Homei (2009) following Johnson and Kotz
(1990) extended Van Assche's results. They put $X_1,\cdots,X_n$ to
be independent, and considered
$$ S_n = R_{1}X_1 + R_2X_2 + \cdots
+R_{n-1}X_{n-1} + R_nX_n,\;\;\;\; n\geq 2\;,  $$ where random
proportions are $R_{i}=U_{(i)}-U_{(i-1)},\; i=1,...,n-1$ and $R_n=
1-\sum_{i=1}^{n-1} R_i$,  $U_{(1)},...,U_{(n-1)}$ are order
statistics from a uniform distribution on $[0,1]$, and $U_{(0)}=0$.
These random proportions are uniformly distributed over the unit
simplex. They employed Stieltjes transform and derived that: (i)
$S_{n}$ possesses the same distribution as $X_{1}$,...,$X_{n}$ if
and only if $X_{1}$,...,$X_{n}$ are degenerated or have a Cauchy
distribution; and (ii) Van Assche's (1987) result for Arcsin holds
for $Z$ only.

In this paper, we introduce two families of distributions, suggested
by an anonymous referee of the article, to whom the author expresses
his deepest gratitude. We say that $Z_1$ is a random variable
between two independent random variables with power distribution, if
the conditionally distribution of $Z_1$ given at $X_1=x_1, X_2=x_2$
is
$$
F_{Z_{1}|x_1,x_2}(z)=\left\{
\begin{array}{cc}
  1  & z\geq {\rm max}(x_1,x_2), \\
  (\frac{z-x_1}{x_2-x_1})^n  & x_1<z<x_2, \\
  1-(\frac{z-x_1}{x_2-x_1})^n  & x_2<z<x_1, \\
  0 & z\leq {\rm {min}}(x_1,x_2).
\end{array}
\right.\eqno(1.1)
$$
The distribution $F_{Z_1|x_1,x_2}$ will be said to follow a
conditionally directed power distribution, when $n$ is an integer.
For $n=1$, the distribution given by (1.1) simplifies to the
distribution $Z$ that was introduced by Van Assche (1987). For
$n=2$, we call $Z_1$ directed triangular random variable. For
further generalizing Van Assche results, we introduce a seemingly
more natural conditionally power distribution.
 We call $Z_2$  two-sided power (TSP) random variable, if the
conditionally distribution of $Z_2$ given at $X_1=x_1, X_2=x_2$ is
$$
F_{Z_{2}|x_1,x_2}(z)=\left\{
\begin{array}{cc}
  1  & z\geq y_2, \\
  (\frac{z-y_1}{y_2-y_1})^n  & y_1<z<y_2, \\
  0 & z\leq y_1.
\end{array}
\right.\eqno(1.2)
$$
The distribution $F_{Z_2|x_1,x_2}$ will be said to follow a
conditionally undirected power distribution, when $y_{1}={\rm
min}(x_1,x_2), y_{2}={\rm max}(x_1,x_2)$ and $n$ is an integer.  For
$n=2$, we call $Z_2$ undirected triangular random variable.\\
Again, for $n=1$, the distribution given by (1.1) simplifies to the
distribution $Z$ that was introduced by Van Assche (1987). The main
aim of this article is providing a couple of generalizations to the
results of Van Assche (1987) for some other values of $n$ (other
than $n=1$). This article is organized as follows. We introduce
preliminaries and previous works in section 2. In section 3, we give
some characterizations for distribution $Z_1$ given in (1.1), when
$n=2$. In section 4, we find distribution of $Z_2$ given in (1.2) by
direct method, and give some examples of such distributions.

\section{Preliminaries and previous works}
In this section, we first review some results of Van Assche (1987)
and then modify them a little bit to fit in our framework, to be
introduced in the forthcoming sections.

Using the Heaviside function ($U(x)=0,\;x<0,\; =1,\; x\geq0$) we
conclude that for any given distinct values $x_{1}$ and $x_{2} $,
the conditional distribution $F_{Z_1|x_1,x_2}(z)$ in (1.1) is
$$F_{Z_1|x_1,x_2}(z)=(\frac{z-x_1}{x_2-x_1})^{n}U(z-x_1)-\sum_{i=1}^{n}{n\choose
i}(\frac{z-x_2}{x_2-x_1})^iU(z-x_2). \eqno(2.1)$$

\textbf{Lemma 2.1.} For distinct reals $x_{1},x_{2},z$ and integer
$n$, we have
$$\frac{-1}{(z-x_1)(x_2-x_1)^n}+\frac{(-1)^n}{(n-1)!}\frac{d^{n-1}}{dx_{2}^{n-1}}(\frac{1}{z-x_2}.\frac{1}{(x_1-x_2)})=\frac{1}{(x_1-z)(x_2-z)^n}.$$\\

\textbf{Proof.} It easily follows from the Leibniz formula.
\hfill$\Box$

\noindent Another tool for proving our main theorem is the following
formula taken from the Schwartz distribution theory, namely,
$$\int_{-\infty}^{\infty}\varphi(x)\Lambda^{[n]}(dx)=\frac{(-1)^{n}}{n!}
\int_{-\infty}^{\infty}\frac{d^{n}}{dx^{n}}\varphi(x)\Lambda(dx),
 \eqno(2.2)$$
where $\Lambda$ is a distribution function and $\Lambda^{[n]}$ is
the $n$-th distributional  derivative of $\Lambda$.\\ The
conditional distribution $F_{Z_1|x_{1},x_{2}}(z)$ given by (1.1)
leads us to a linear functional on complex-valued functions
$f:\mathbb{R}\rightarrow \mathbb{C}$, defined on the set of real
numbers $\mathbb{R}$:
$$F_{Z_{1}|x_{1},x_{2}}(f)=\frac{f(x_1)}{(x_2-x_1)^n}-\sum_{i=1}^{n}\frac{1}{(n-i)!(x_2-x_1)^{i}}\frac{d^{n-i}}{dz^{n-i}}f(x_2).$$
It easily follows that
$$F_{Z_{1}|x_{1},x_{2}}(af+bg)=aF_{Z_{1}|x_{1},x_{2}}(f)+bF_{Z_{1}|x_{1},x_{2}}(g), \eqno(2.3) $$
for any complex-valued functions $f,g$ and  complex constants $a,b$.
We note that $F_{Z_{1}|x_{1},x_{2}}(z)=F_{Z_{1}|x_{1},x_{2}}(f_z)$,
whenever $f_{z}(x)=(z-x)^{n}U(z-x)$ and\\
$$F_{Z_{1}|x_{1},x_{2}}(f_z)=\frac{f_{z}(x_1)}{(x_2-x_1)^n}-\sum_{i=1}^{n}\frac{1}{(n-i)!(x_2-x_1)^{i}}\frac{d^{n-i}}{dz^{n-i}}f_{z}(x_2).$$
Also we note that $U(z-x)= \frac{(-1)^{n}}{(n)!}
\frac{d^{n}}{dx^{n}} f_{z}(x).$ Thus
$$P(Z_{1}\leq z)=\int_{\mathbb{R}}U(z-x)dF_{Z_{1}}(x)=
\int_{\mathbb{R}^{2}}F_{Z_{1}|x_{1},x_{2}}(z)\prod_{i=1}^{2}
F_{X_{i}}(dx_{i}),$$
 can be viewed as:
$$\int_{\mathbb{R}}\frac{(-1)^{n}}{(n)!}
\frac{d^{n}}{dx^{n}}f_{z}(x)dF_{Z_1}(x)=\int_{\mathbb{R}^{2}}F_{Z_{1}|x_{1},x_{2}}(f_z)
\prod_{i=1}^{2}F_{X_{i}}(dx_{i}).\eqno(2.4)
$$ Therefore by using (2.3) along with (2.4) and a standard
argument in the integration theory, we obtain that
$$\int_{\mathbb{R}}\frac{(-1)^{n}}{(n)!}
\frac{d^{n}}{dx^{n}}f(x)dF_{Z_1}(x)=\int_{\mathbb{R}^{2}}F_{Z_{1}|x_{1},x_{2}}(f)
\prod_{i=1}^{2}F_{X_{i}}(dx_{i}), \eqno(2.5)$$ for any infinitely
differentiable functions $f$ for which the corresponding integrals
are finite. Now (2.5) together with (2.2) lead us to
$$\int_{\mathbb{R}}f(x)dF_{Z_1}^{(n)}(x)=\int_{\mathbb{R}^{2}}F_{Z_{1}|x_{1},x_{2}}(f)
\prod_{i=1}^{2}F_{X_{i}}(dx_{i}), \eqno(2.6)$$
for the above mentioned functions $f$, where $F_{Z_1}^{(n)}$ is the $(n)$-th distributional derivative of the distribution of $Z_1$.\\
 Let us denote the
Stieltjes transform of a distribution $H$ by
$${\cal S}(H,z)=\int_{\mathbb{R}}\frac{1}{z-x}H(dx),$$
for every $z$ in the set of complex numbers $\mathbb{C}$ which
 does not belong to the support of $H$, i.e., $z \in \mathbb{C}\cap
 (\mbox{supp}H)^{\cal C}$. For more on the Stieltjes transform, see Zayed
 (1996).

The following lemma indicates how the Stieltjes transform of $Z_1$,
and $X_1,X_2$ are related.

\textbf{Lemma 2.2. } Let $Z_1$ be a random variables that satisfies
(1.1). Suppose that the random variables $X_1$ and $X_2$ are
independent and continuous with distribution functions $F_{X_1}$ and
$F_{X_2}$, respectively. Then
$$\frac{1}{n}{\cal S}^{(n)}(F_{Z_1},z)=-{\cal S}(F_{X_1},z){\cal
S}^{(n-1)}(F_{X_2},z), \;\;\;z\in \mathbb{C}\bigcap_{i=1}^2
(\mbox{supp}F_{X_i})^{\cal C}.$$

\textbf{Proof.} It follows from (2.6) that
$${\cal S}(F_{Z_{1}}^{(n)},z)=\int_{\mathbb{R}^{2}}
F_{Z_{1}|x_{1},x_{2}}(g_z) \prod_{i=1}^{2}F_{X_{i}}(dx_{i}),$$ and
$$\frac{1}{n!}\frac{d^{n}}{dz^{n}}{\cal
S}(F_{Z_{1}},z)=\int_{\mathbb{R}^{2}}F_{Z_{1}|x_{1},x_{2}}(g_z)\prod_{i=1}^{2}
F_{X_{i}}(dx_{i}),$$
 for $g_{z}(x)=\frac{1}{z-x}$. Now,
it follows that
$$F_{Z_{1}|x_{1},x_{2}}(g_z)=\frac{\frac{1}{z-x_1}}{(x_2-x_1)^n}-\sum_{i=1}^{n}\frac{1}{(n-i)!(x_2-x_1)^i}\frac{d^{n-i}}{dz^{n-i}}\frac{1}{z-x_2},$$
and by using Lemma 2.1, we have
\begin{eqnarray*}
F_{Z_{1}|x_{1},x_{2}}(g_z)&=&\frac{(-1)^n}{(z-x_1)(z-x_2)^n}.
\end{eqnarray*}
Therefore,
$$\frac{1}{n!}\frac{d^{n}}{dz^{n}}{\cal S}(F_{Z_1},z)=\int_{\mathbb{R}^{2}}
\frac{(-1)^n}{(z-x_1)(z-x_2)^n} \prod_{i=1}^{2} F_{X_{i}}(dx_{i}),
$$ and
$$\frac{1}{n}{\cal S}^{(n)}(F_{Z_1},z )=-{\cal S}(F_{X_{1}},z){\cal S}^{(n-1)}(F_{X_{2}},z), \;\;\;z\in \mathbb{C}\bigcap_{i=1}^2
(\mbox{supp}F_{X_i})^{\cal C}.\eqno(2.7)$$ This finishes the proof.
 \hfill$\Box$

Note that Van Assche's lemma is the case of $n=1$:
$$-{\cal S}^{'}(F_{Z_1},z)={\cal S}(F_{X_{1}},z){\cal S}(F_{X_{2}},z).$$
We also note that the Stieltjes transform of Cauchy distribution,
i.e., ${\cal S}(F,z)=\frac{1}{z+c}$, satisfies (2.7).

\section{Directed triangular random variable}
Let us now review Van Assche's result for directed triangular random
variables.

\textbf{Theorem 3.1.} If $X_1$ and $X_2$ are independent random
variables with a common distribution $F_X$, then the
characterizations of $Z_1$ for $n=1$ and $n=2$ are identical.

\textbf{Proof.} We note that $X_1$ and $X_2$ have a common
distribution function $F_X$. By using Lemma 2.2 for $n=2$, we have
$$-\frac{1}{2}{\cal S}^{''}(F_{Z_1},z)={\cal S}(F_X,z){\cal
S}^{'}(F_X,z),$$ and so $$-{\cal
S}^{''}(F_{Z_1},z)=\frac{d}{dz}{\cal S}^{2}(F_X,z),$$ and $$-{\cal
S}^{'}(F_{Z_1},z)={\cal S}^{2}(F_X,z). \eqno(3.1)$$ We note that the
Stieltjes transform tends to zero, when $z$ is sufficiently large.
In that case the constant in the differential equation will be zero.
The equation (3.1) is exactly the equation obtained by Van Assche
(1987) when $X_1$ and $X_2$ have a common distribution; so his
results hold in our framework as well. \hfill$\Box$ \\ This clever
proof is due to the anonymous referee. Now, we apply Lemma 2.2 for
some characterizations, when $X_1$ and $X_2$ are not identically
distributed.

\textbf{Theorem 3.2.} Let $X_1$ and $X_2$ be independent random
variables and $Z_1$ be a directed triangular random variable
satisfying $(1.1)$. For $n=2$, we have,
 \vskip 3mm

(a) if $X_1$ has uniform distribution on $[-1,1]$,
 then $Z_1$ has semicircle distribution on $[-1,1]$ if and only if $X_2$ has Arcsin distribution on
 $[-1,1]$;

(b) if $X_1$ has uniform distribution on $[-1,1]$, then $Z_1$ has
power semicircle distribution if and only if $X_2$ has power
semicircle distribution, i.e.,
$$f(z)=\frac{3(1-z^2)}{4}\;,\;\;\;\;-1\leq z\leq 1;$$

(c) if $X_1$ has Beta$(1,1)$ distribution on $[0,1]$, then $Z_1$ has
Beta$(\frac{3}{2},\frac{3}{2})$ distribution if and only if $X_2$
has Beta$(\frac{1}{2},\frac{1}{2})$ distribution;

(d) if $X_1$ has uniform distribution on $[0,1]$, then $Z_1$ has
Beta$(2,2)$ distribution if and only if $X_2$ has Beta$(2,2)$
distribution.

\textbf{Proof.} (a) For the ``if" part we note that the random
variable $X_1$ has uniform distribution and $X_2$ has arcsin
distribution on [-1,1]; then
$${\cal S}(F_{X_1},z)=\frac{1}{2}(\ln|z+1|-\ln|z-1|),$$ and
$${\cal S}(F_{X_2},z)=\frac{1}{\sqrt{z^2-1}}.$$
From Lemma 2.2 and substituting the corresponding Stieltjes
transforms of distributions, we get
$${ \cal S}^{''}(F_{Z_1},z)=\frac{2}{(z^2-1)^{\frac{3}{2}}}.$$
The solution ${\cal S}(F_{Z_1},z)$ is
$${\cal S}(F_{Z_1},z)=2(z-\sqrt{z^2-1}),$$
which is the Stieltjes transform of the semicircle distribution on
$[-1,1]$.\\ For the ``only if" part we assume that the random
variable $Z_1$ has semicircle distribution. Then it follows from
lemma 2.2 that
$${\cal S}(F_{X_2},z)\; \frac{1}{1-z^2} = \frac{-1}{(z^2-1)^{\frac{3}{2}}}.$$
 The proof is completed.

(b) By an argument similar to that given in (a) and solving the
following differential equations,
$$S^{''}(F_{Z},z)=\frac{2}{(z^2-1)}(\frac{3z}{2}+\frac{3}{4}(1-z^2)({\rm ln}|z+1|-{\rm ln}|z-1|)), \ \textrm{(for
 the  ``if"   part), and}$$
$$\frac{1}{1-z^2}S(F_{X_2},z)=\frac{3}{4}\frac{2z+(1-z^2)({\rm ln}|z+1|-{\rm ln}|z-1|)}{(1-z^2)}, \ \textrm{(for the
``only if" part)},$$ the proof can be completed.

(c) By Lemma (2.2), we have
$$-\frac{1}{2}S^{''}(F_{Z},z)=\frac{-1}{z(z-1)}\frac{1}{\sqrt{z(z-1)}}, \ \textrm{(for
 the  ``if"   part), and}$$
$$\frac{-1}{z(z-1)}S(F_{X_2},z)=\frac{-1}{z(z-1)\sqrt{z(z-1)}}, \ \textrm{(for the
``only if" part)}.$$ The proof can be completed by solving the above
differential equations.

(d) By Lemma (2.2), we have
$$\hfill{\cal S}^{''}(F_{Z_1},z)= \frac{-2}{z(z-1)}(6(z^2-z)(\ln|z|-\ln|z-1|)-6z+3), \ \textrm{(for
 the  ``if"   part), and}$$
$${\cal S}(F_{X_2},z) = 6(z-z^2)(\ln|z|-\ln|z-1|)+6z-3, \ \textrm{(for the
``only if" part)}.$$ Solving the differential equations, can
complete the proof. \hfill$\Box$

\section{TSP random variables}
In section 3, we used a powerful method, based on the use of
Stieltjes transforms, to obtain the distribution of $Z_1$ given in
(1.1). It seems that one can not use that method to find the
distribution of $Z_2$ given in (1.2). So we employ a direct method
to find the distribution of $Z_2$.

Let us follow Lemma 4.1 to find a simple method to get the
distribution of $Z_2$. The work of Soltani and Homei (2009b) leads
us to the following lemma.

\textbf{Lemma 4.1.} Suppose $W$ has a power distribution with
parameter $n$, $n\geq 1$, $n$ is an integer, and let
$Y_1={\rm{Min}}(X_1,X_2)$, $Y_2={\rm Max}(X_1,X_2),$ where $X_1$ and
$X_2$ are independent random variables. Let
$$X=Y_1+W(Y_2-Y_1). \eqno(4.1)$$
Then

(a) $X$ is a TSP random variable.

(b) $X$ can be equivalently defined by
$$X=\frac{1}{2}(X_1+X_2)+(W-\frac{1}{2})|X_1-X_2|.$$

\textbf{Proof.} (a)\begin{eqnarray*}
F_{X|x_1,x_2}(z)&=& P(Y_1+W(Y_2-Y_1)\leq z|X_1=x_1,X_2=x_2)\\
&=&P(y_1+W(y_2-y_1)\leq z)\\
&=&(\frac{z-y_1}{y_2-y_1})^n.
\end{eqnarray*}

(b) The proof can be completed by substuting   ${\rm Min}(X_1,X_2)$
and ${\rm Max}(X_1,X_2)$ with $Y_1$ and $Y_2$ in (4.1). \hfill$\Box$

\subsection{Moments of TSP random variables}
The following theorem provides equivalent conditions for
$\mu_{k}^{'}=EZ_{2}^{k}$.

\textbf{Theorem 4.1.1.} Suppose that $Z_2$ is a TSP random variable
satisfing (1.2). If $X_1$
and $X_2$ are random variables and $E|X_i|^{k}<\infty$, $i=1,2$, for all integers $k$, then\\
(a) $EZ_{2}^{k}=n\frac{\Gamma(k+1)}{\Gamma(k+n+1)}\sum_{i=0}^{k}\frac{\Gamma(k-i+n)}{\Gamma(k-i+1)}E(Y_{1}^{i}Y_{2}^{k-i})$;\\
(b) $EZ_{2}^{k}=\sum_{i=0}^{k}{k\choose
i}(\frac{1}{2})^{k-i}E(W-\frac{1}{2})^{i}E(X_1+X_2)^{k-i}|X_1-X_2|^{i}$;\\
(c) $EZ_{2}^{k}=\sum_{i=0}^{k}{k\choose
i}\frac{n}{n+i}E(Y_{1}^{k-i}(Y_2-Y_1)^i)$.

\textbf{Proof.} (a) By using Lemma 4.1, we obtain that
\begin{eqnarray*}
EZ_{2}^{k}&=& E(\sum_{i=0}^{k}{k \choose i}(1-W)^{i}Y_{1}^{i}W^{k-i}Y_{2}^{k-i})\\
&=&\sum_{i=0}^{k}{k \choose i}E(W^{k-i}(1-W)^{i})E(Y_{1}^{i}Y_{2}^{k-i})\\
&=&n\frac{\Gamma(k+1)}{\Gamma(k+n+1)}\sum_{i=0}^{k}\frac{\Gamma(k-i+n)}{\Gamma(k-i+1)}EY_{1}^{i}Y_{2}^{k-i}.\\
\end{eqnarray*}

(b) This can be easily proved by Lemma 4.1(b).

(c) It straightforwardly follows from (4.1). \hfill$\Box$

Let us consider expectation and variance of $Z_2$. First, we suppose
that ${\rm E}Y_1=\mu_1$, ${\rm E}Y_2=\mu_2$, ${\rm Var}
Y_1=\sigma_1^{2}$, ${\rm Var} Y_2=\sigma_{2}^{2}$, and ${\rm
Cov}(Y_1,Y_2)=\sigma_{12}$. Then
$$EZ_2=\frac{\mu_1+n\mu_2}{n+1},$$
and also, if ${\rm E}X_1={\rm E}X_2=0$, then
$$E(Z_2)=EY_{1}+\frac{n}{n+1}(EY_{2}-EY_{1}).$$
By $Y_{1}+Y_{2}=X_{1}+X_{2}$, we have
$$E(Z_2)=E(Y_1)+\frac{n}{n+1}(-2EY_{1})=\frac{1-n}{1+n}EY_1.\eqno(4.2)$$
It can easily follow from (4.2) that the Arcsin result of Van Assche
(1987) is only true for $n=1$, and also, one can see that Theorem
(3.2) in section 3 does not hold for the above $Z_2$.\\
About the variance, we have $$Var
Z_2=\frac{n(\mu_{1}-\mu_2)^2+n(n+1)^2
\sigma_{2}^{2}+2(n+1)(\sigma_{1}^{2}+n
\sigma_{12})}{(n+1)^2(n+2)}.$$ Following the computation of
expectation and variance, we evaluate them for some well-known
distributions. If $X_1$ and $X_2$ have standard normal
distributions, then from Theorem 4.1.19b) and the fact that
$X_1-X_2$ and $X_1+X_2$ are independent, it follows that their
first, second and third order moments are equal, respectively, to
$$EZ_2=\frac{1}{\sqrt{\pi}}(\frac{n-1}{n+1}),$$
$$EZ_{2}^{2}=\frac{n^2+n+2}{(n+1)(n+2)}, \;\; {\rm and}$$
$$EZ_{2}^{3}=\frac{1}{2\sqrt{\pi}}\frac{5n^3+12n^2+13n-30}{(n+3)(n+2)(n+1)}.$$
Also, in case $X_1$ and $X_2$ have uniform distributions, Theorem
4.1.1(b) implies that,
$$EZ_{2}^{k}=n\frac{\Gamma(k+1)}{\Gamma(n+k+1)}\sum_{i=0}^{k}\frac{\Gamma(k-i+n)}{\Gamma(k-i+1)}\frac{2}{(k+2)(i+1)},$$
$$EZ_{2}=\frac{2n+1}{3(n+1)}, \;\; {\rm and}$$
$$Var(Z_2)={\frac{1}{18}}\frac{n^3+3n^2+6n+2}{(n+1)^2(n+2)}.$$
Since some distributions do not have any moments, Theorem 4.1.1 is
not applicable for investigating Van Assche's results for them,
whence, we prove the following theorem:

\textbf{Theorem 4.1.2.} Suppose that $Z_2$ is a TSP random variable
satisfying (4.1). Then\\
(a) $Z_2$ is location invariant;\\
(b) if $X_1$ and $X_2$ have symmetric distribution around $\mu$,
then $Z_2$ has symmetric distribution around $\mu$, only when $n=1$.

\textbf{Proof.}

(a) Is immediate.

(b) We can assume without loss of generality that $\mu=0$. If $Z_2$
has a symmetric distribution around zero, then
$$Y_{1}+W(Y_{2}-Y_{1})\stackrel{d}{=}-[Y_{1}+W(Y_2-Y_1)].$$
We note that
$$Y_{1}+W(Y_{2}-Y_{1})\stackrel{d}{=}[-Y_{1}+W(-Y_2-(-Y_1))].$$
Since, $-{\rm {Min}}(X_1,X_2)={\rm {Max}}(-X_1,-X_2)$,
$X_1\stackrel{d}{=}-X_1$ and $X_2\stackrel{d}{=}-X_2$, we have
$$Y_{1}+W(Y_{2}-Y_{1})\stackrel{d}{=}Y_{2}+W(Y_1-Y_2).\eqno(4.3)$$
By equating the conditional distributions given at $X_1=x_1$ and
$X_2=x_2$ in (4.2), we conclude that $n=1$. \hfill$\Box$\\
It can also easily follow from Theorem (4.1.1) that the Cauchy
result of Van Assche (1987) is true only for $n=1$.

\subsection{Distributions of TSP random variables}
In this subsection, we investigate computing distributions by the
direct method. We will give two examples of derivation based on
(4.1). This method may be complicated in some cases, but we have
chosen some easy to follow examples.

\textbf{Example 4.2.1.} Let $X_1,X_2$ and $W$ be independent random
variables such that $X_1$ and $X_2$ are uniformly distributed over
$[0,1]$, and $W$ has a power function distribution with parameter
$n$. We find the value $f_{Z_2}(z;n)$ by means of $f_{Z_2|W}(z|w)$;
therefore
$$
f_{Z_2|W}(z|w)=\left\{
\begin{array}{cc}
\frac{2z}{w},\;\;0<z<w,\\
\frac{2(1-z)}{1-w},\;\;w<z<1.\\
\end{array}
\right.\eqno(4.4)
$$
By using the distribution of $W$, the density function
$f_{Z_2}(z;n)$ can be expressed in terms of the Gauss hypergeometric
function $F(a,b,c;z)$, which is a well-known special function.
Indeed according to Euler's formula, the Gauss hypergeometric
function assumes the integral representation
$$F(a,b,c;z)=\frac{\Gamma(c)}{\Gamma(b)\Gamma(c-b)}\int_{0}^{1}t^{b-1}(1-t)^{c-b-1}(1-tz)^{-a}dt,$$
where $a,b,c$ are parameters subject to $-\infty<a<+\infty$,
$c>b>0$, whenever they are real, and $z$ is the variable (see Zayed
1996). By using Euler's formula, the density function of $Z_2$ can
be expressed as follows:
$$f_{Z_2}(z;n)=\frac{2nz}{n-1}(1-z^{n-1})+2(1-z)z^{n}F(1,n,n+1,z),\;0<z<1,\eqno(4.5)$$
where $n>0$ and $n\neq 1$. When $n=1$, similar calculations lead to
the following distribution
$$f_{Z_2}(z)=-2(1-z){\rm log}(1-z)-2z {\rm log}(z),\;\;0<z<1.$$

The probability density function $f_{Z_2}(z)$ was introduced by
Johnson and Kotz (1990), for the first time, under the title
``uniformly randomly modified tine". So $f_{Z_2}(z;n)$ can be seen
as an extension of the above mentioned distribution. We note that,
from (4.1) and a simple Monte Carlo procedure using only simulated
uniform variables, one is able to simulate the distribution (4.5).

\textbf{Example 4.3.1.} Let $X_1$ and $X_2$ be independent random
variables with Beta$(1,2)$ distribution. Then if $W$ has Beta$(3,1)$
distribution, $Z_1$ has Beta$(2,3)$ distribution.

In the following theorem we compute the Stieltjes transform of $Z_2$
for $n=2$. Let us remark that the complexity of the integral in the
theorem indicates that for this case the direct method is preferred.

\textbf{Theorem 4.4.1} Let $Z_2$ be a undirected triangular random
variable that satisfies (1.2). Suppose that the random variables
$X_1$ and $X_2$ are independent and continuous with the distribution
functions $F_{X_1}$ and $F_{X_2}$, respectively. Then
$$-\frac{1}{2}{\cal S}^{'''}(F_Z,z)={\cal S}^{'}(F_{X_1},z){\cal S}^{'}(F_{X_2},z)+2{\cal S}(F_{X_1},F_{X_2},z),$$
where
$${\cal S}(F_{X_1},F_{X_2},z)=\int_{\mathbb{R}^{2}}\frac{1}{(z-x_1)(z-x_2)(x_2-x_1)^2}\prod_{i=1}^{2} F_{X_i}(dx_i).$$

\textbf{Proof.} By using an argument similar to that given in
Section 3, we can conclude that
$$\int f(x)dF_{Z_2}^{(2)}(x)=\int_{\mathbb{R}^{2}}
F_{Z_{2}|x_{1},x_{2}}(f) \prod_{i=1}^{2}F_{X_{i}}(dx_{i}).$$ So,
$$-\frac{1}{2}{\cal S}^{'''}(F_{Z_2},z)=\int_{\mathbb{R}^{2}}
F_{Z_{2}|x_{1},x_{2}}(g_z) \prod_{i=1}^{2}F_{X_{i}}(dx_{i}),$$ for
$g_{z}(x)=\frac{1}{(z-x)^2}$. From
$$F_{Z_2|x_1,x_2}(g_z)=\frac{\frac{1}{(z-x_1)^2}}{(x_2-x_1)^2}+\frac{\frac{1}{(z-x_2)^2}}{(x_1-x_2)^2}$$
and by using partial fractional rule, we have
$$F_{Z_2|x_1,x_2}(g_z)=\frac{1}{(z-x_1)^{2}(z-x_2)^2}+\frac{2}{(x_2-x_1)^2}\frac{1}{(z-x_1)(z-x_2)}.$$
Therefore,
$$-\frac{1}{2}{\cal S}^{'''}(F_{Z_2},z)=\int_{\mathbb{R}^{2}}(\frac{1}{(z-x_1)^{2}(z-x_2)^2}+\frac{2}{(x_2-x_1)^{2}(z-x_1)(z-x_2)})\prod_{i=1}^{2}
F_{X_{i}}(dx_{i}),$$ and
$$-\frac{1}{2}{\cal S}^{'''}(F_{Z_2},z)={\cal S}^{'}(F_{X_1},z){\cal S}^{'}(F_{X_2},z)+2{\cal S}(F_{X_1},F_{X_2},z).$$
This finishes the proof.
 \hfill$\Box$

It is worth mentioning that the present method yields other
extensions too; the following is such an example.

\textbf{Example 4.3.2.} Suppose that $X_1,X_2,W$ are independent
random variables. If $X_1$ and $X_2$ have uniform distributions on
$[0,1]$ and $W$ has Beta$(2,2)$ distribution, then $Z_2$ has the
same distribution as $W$.

If the product moments of order statistics are known, those of $W$
can be derived from that of $Z_2$ by using Theorem 4.1.1(a). Then
the distribution of $W$ is characterized by that of $Z_2$.

By an argument similar to the one given in Example 4.2.1, when $W$
has a Beta distribution with parameters $n$ and $m$, we find the
distribution $f_{Z_2}(z;n,m)$ as
$$\frac{B(n-1,m)}{B(n,m)}2z(1-I_{z}(n-1,m))+\frac{B(n,m-1)}{B(n,m)}2(1-z)I_{z}(n,m-1),\;\;0<z<1,$$
where $I_{x}(a,b)$ is incomplete Beta function:
$$I_{x}(a,b)=\frac{1}{B(a,b)}\int_{0}^{x}t^{a-1}(1-t)^{b-1}dt,\;\;  (a,b>0).$$

\section{Conclusions}
We have described how (a) methods of Stieltjes transform, and (b)
directed methods, could be used for obtaining the distributions,
characterizations and properties of the random mixture of variables
defined in (1.1) and (1.2). Of course each one of the methods (a) or
(b) has its own advantages and disadvantages, and none of them has a
preference over the other. The TSP random variable when $X_1$ and
$X_2$ have uniform distributions, led us to a new family of
distributions which can be regarded as some generalization of
``uniformly randomly modified tine". The proposed model in the
direct method can easily lead to distribution generalizations,
though this is not possible for the first method, but here the
characteristics can be easily computed.

\section{Acknowledgment}
The author is deeply grateful to the anonymous referee for reading
the original manuscript very carefully and for making valuable
suggestions.

\end{document}